\documentclass[11pt]{amsart}

\usepackage{amsmath}
\usepackage{amssymb}
\usepackage{graphicx}
\usepackage{xcolor}
\newtheorem{theorem}{Theorem}[section]
\newtheorem{corollary}[theorem]{Corollary}
\newtheorem{lemma}[theorem]{Lemma}
\newtheorem{proposition}[theorem]{Proposition}
\theoremstyle{definition}

\newtheorem{remark}[theorem]{Remark}

\numberwithin{equation}{section}

\title[ On the entropy numbers and  the Kolmogorov widths]{On the entropy numbers and  the Kolmogorov widths}

\author[G. Petrova]{Guergana Petrova}
\address[G. Petrova]{Department of Mathematics, Texas A$\&$M University, College Station, TX 77843, USA}
\email{{\tt gpetrova$@$math.tamu.edu}}

\author[P. Wojtaszczyk]{Przemys{\l}aw Wojtaszczyk}
\address[P. Wojtaszczyk]{Institut of Mathematics, Polish Academy of Sciences, ul. {\'S}niadeckich 8,  00-656 Warszawa, Poland}
\email{\tt wojtaszczyk$@$impan.pl}

\keywords{Kolmogorov widths, entropy numbers\\
  Dedicated to Ron DeVore, with the utmost respect and admiration}


{

\begin{document}

\begin{abstract}
Direct estimates between linear or nonlinear Kolmogorov widths and entropy numbers are presented. These
estimates are derived using the recently introduced Lipschitz widths. { Applications for $m$-term approximation are obtained.}
\end{abstract}

\maketitle


\section{Introduction}

We consider a Banach space $(X,\|\cdot\|_X)$ (or a Hilbert space $H$) equipped with a norm $\|\cdot\|_X$  and a compact subset ${\mathcal K}\subset X$  of $X$. Typically, ${\mathcal K}$ is 
a  finite ball in smoothness spaces like the Lipschitz, Sobolev, or Besov spaces. 

A well known classical result, called the Carl's inequality, see     \cite{C} or  \cite{LGM}, compares a certain
characteristic of   the set ${\mathcal K}$, called  entropy numbers $e_k({\mathcal K})_X$, with  its approximability by linear spaces, measured by its Kolmogorov width $d_{k}({\mathcal K})_X$. The Carl's  inequality  states that 
for each $r>0$, there is a constant $C(r)$ such that for all $n\in \mathbb N$,
\begin{equation}
 \label{carl}
\max_{1\leq k\leq n}k^re_k({\mathcal K})_X\leq C(r)\max_{1\leq m\leq n}m^rd_{m-1}({\mathcal K})_X.
\end{equation}
Inequality (\ref{carl}) has been generalized in \cite{T}, where the nonlinear Kolmogorov widths $d_n({\mathcal K},N)_X$ have been used instead of the linear Kolmogorov widths $d_{k}({\mathcal K})_X$.  More precisely, it has been shown there that for each $r>0$, there is a constant $C(r,\lambda)$ such that for all $n\in {\mathbb N}$, 
\begin{equation}
 \label{carl1}
\max_{1\leq k\leq n}k^re_k({\mathcal K})_X\leq C(r,\lambda)\max_{1\leq m\leq n}m^rd_{m-1}({\mathcal K},\lambda^m)_X,
\end{equation}
with $\lambda>1$ a fixed constant.
In addition, it was also proven that for each $r>0$, there is a constant $C(r,a)$ such that for all $n\in {\mathbb N}$, 
\begin{equation}
 \label{carl11}
\max_{1\leq k\leq n}k^re_{(a+r)k\log k}({\mathcal K})_X\leq C(r,a)\max_{1\leq m\leq n}m^rd_{m-1}({\mathcal K},m^{am})_X, 
\end{equation}
where $a>0$ is a fixed constant and  $k\log k$ cannot be replaced by a slower growing function of $k$.

All these inequalities are primarily useful  when  the linear or nonlinear 
Kolmogorov widths  decay as a power of $m$. In this paper,  we give  finer extensions of the (generalized) Carl's inequalities 
(\ref{carl}), (\ref{carl1}) and (\ref{carl11}), using  the recently introduced in \cite{PW} Lipschitz widths. We start with some definitions, presented in \S
 \ref{S2}, and continue, see  \S\ref{S3}, with a comparison between  the nonlinear Kolmogorov widths and the Lipschitz widths. Our main results are presented in \S\ref{S4}, where we give a direct comparison between the entropy numbers of ${\mathcal K}$ 
 and its linear and nonlinear Kolmogorov widths. 
  Finally, in  \S \ref{appl}, we derive what these estimates mean for the $m$-term approximation in Hilbert spaces.

\section{Preliminaries}
\label{S2}
  We start this section with the definition of  Kolmogorov widths.
  If we fix the value of $n\ge 0$, the {\it  Kolmogorov $n$-width}  $d_n({\mathcal K})_X$ of ${\mathcal K}$
 is defined as
 \begin{equation}
\nonumber
  d_0({\mathcal K})_X:=\sup_{f\in {\mathcal K}} \|f\|_X, \quad d_n({\mathcal K})_X:=\inf_{\dim(X_n)=n} \sup_{f\in {\mathcal K}}{\rm dist}(f,X_n) _X, \quad  n\geq 1,
 \end{equation}
 where the infimum is taken over all linear spaces $X_n\subset X$ of dimension $n$.
  These are the  classical Kolmogorov widths introduced in \cite{K}, or consult \cite{LGM}
  for their modern exposition. To distinguish them from the introduced later nonlinear Kolmogorov widths,   we call them {\it linear Kolmogorov $n$-widths}.
They describe  the optimal performance possible  for the approximation
 of the model class ${\mathcal K}$ using linear spaces of dimension $n$. However, they do not tell us how to select
 a (near) optimal space $Y$ of dimension $n$ for this purpose. 
Let us also note that in the definition of Kolmogorov width,  we are not requiring that   the mapping 
 which sends $f\in {\mathcal K}$ into an approximation to $f$  is a linear map.    
 
 A generalization of this concept was introduced in \cite{T}, where the so called {\it nonlinear Kolmogorov $(n,N)$-width} 
 $d_n({\mathcal K},N)_X$ was defined for $N\geq 1$ as
 $$
d_0({\mathcal K},N)_X:=\sup_{f\in K} \|f\|_X, 
$$
$$
d_n({\mathcal K},N)_X:= \inf_{\mathcal L_N}\sup_{f\in \mathcal K}\inf_{X_n\in {\mathcal L_N}}{\rm dist}(f,X_n) _X, \quad  n\geq 1,
$$
where the last infimum is over the sets ${\mathcal L}_N$ of at most $N$ linear spaces  $X_n\subset X$  of dimension $n$. Note that here  the choice of the linear subspace $X_n\in {\mathcal L_N}$ from which we choose the best approximation to $f$ depends on $f$. 
 Clearly, 
$d_n({\mathcal K},1)_X=d_n({\mathcal K})_X$, and the bigger the $N$ is, the more flexibility we have to approximate $f$.
These nonlinear Kolmogorov widths are  used in estimating from below the best $m$-term approximation,  see e.g. \cite{HH, T}. The cases considered in \cite{T} are  the cases when $N=\lambda^n$, and $N=n^{an}$, where $\lambda>1$ and $a>0$ are fixed constants, respectively. 
A useful observation that we are going to utilize is that both Kolmogorov widths are homogenous. Namely, if $\mathcal{K} \subset X$ and $t\in \mathbb R$, we have
\begin{equation} \label{hom}
d_n({t \mathcal K},N)_X=|t|d_n({\mathcal K},N)_X\  \mbox{ and  }\  d_n(t\mathcal{K})_X =|t| d_n(\mathcal{K})_X,
\end{equation}
 where ${t \mathcal K}:=\{tf:\,f\in \mathcal K\}$.

In going further, we introduce first the 
{\it minimal  $\epsilon$-covering number} $ N_\epsilon(\mathcal K)$ of a compact set $\mathcal K\subset X$.
A collection 
$\{g_1,\ldots,g_{m}\}\subset X$
of  elements of $X$
is called an $\epsilon$-covering  of $\mathcal K$ if 
$$
\mathcal K\subset \bigcup_{j=1}^m B(g_j,\epsilon), \quad \hbox{where}\quad B(g_j,\epsilon):=\{f\in X:\,\|f-g_j\|_X \leq \epsilon \}.
$$
An  $\epsilon$-covering of $\mathcal K$ whose cardinality   is minimal is called {\it minimal $\epsilon $-covering} of $\mathcal K$. 
We denote by  $N_\epsilon(\mathcal K)$   the  cardinality of the  minimal $\epsilon$-covering  of $\mathcal K$.
{\it  Minimal   inner $\epsilon $-covering number} $\tilde N_\epsilon (\mathcal K)$ of a compact set $\mathcal K\subset X$
is defined exactly as $N_\epsilon(\mathcal K)$ 
but we additionally require that the centers  $\{g_1,\ldots,g_{m}\}$ of the covering are elements from $\mathcal K$. 

{\it Entropy numbers} $e_n({\mathcal K}) _X$, $n\geq 0$, of the compact set ${\mathcal K}\subset X$ are defined 
as the infimum of all $\epsilon>0$ for which $2^n$ balls with centers from $X$ and radius $\epsilon$ cover ${\mathcal K}$.  
If we put the additional restriction that the centers of these balls are from ${\mathcal K}$, then we define the so called {\it   inner entropy numbers} $\tilde e_n({\mathcal K})_X$.
Formally,  we  write
$$ 
e_n({\mathcal K})_X=\inf\{ \epsilon>0 \ :\ {\mathcal K} 
\subset \bigcup_{j=1}^{2^n} B(g_j,\epsilon), \ g_j\in X, \ j=1,\ldots,2^n\},
$$
$$
\tilde e_n({\mathcal K})_X=\inf\{ \epsilon>0 \ :
\ {\mathcal K} \subset \bigcup_{j=1}^{2^n} B(h_j,\epsilon), \ h_j\in {\mathcal K}, \ j=1,\ldots,2^n\}.
$$
\noindent
A collection 
$\{f_1,\ldots,f_{\ell}\}\subset \mathcal K$
of  elements from $\mathcal K$
is called an $\epsilon$-packing of $\mathcal K$ if 
$$
\min_{i\neq j}\|f_i-f_j\|_X>\epsilon.
$$
An  $\epsilon$-packing of $\mathcal K$ whose size   is maximal is called {\it maximal $\epsilon$-packing} of $\mathcal K$.  We denote by $\tilde P_\epsilon(\mathcal K)$ the cardinality of the  maximal $\epsilon$-packing of $\mathcal K$.
We have the following inequalities for every $\epsilon>0$ and every compact set $\mathcal K$
\begin{equation}
\label{known}
\tilde P_\epsilon(\mathcal K)\geq \tilde N_\epsilon(\mathcal K)\geq \tilde P_{2\epsilon}(\mathcal K),
 \end{equation}
 and
\begin{equation}
 \label{innerentropy}
e_n(\mathcal K)_X\leq \tilde e_n(\mathcal K)_X\leq 2 e_n(\mathcal K)_X.
\end{equation}

Finally, we  introduce the {\it Lipschitz widths $d_n^\gamma({\mathcal K})_X$}, $\gamma\geq 0$, $n\geq 1$, of  the 
compact set ${\mathcal K}\subset X$, see \cite{PW}.
We 
denote by    $({\mathbb R}^n,\|.\|_{Y_n})$, $n\geq 1$,
the $n$-dimensional Banach space with a fixed norm $\|\cdot\|_{Y_n}$.
For  $\gamma\geq 0$, we first define 
 the {\it fixed Lipschitz} width $d^\gamma({\mathcal K},Y_n)_X$,
\begin{equation}
\nonumber
d^\gamma({\mathcal K},Y_n)_X:= \inf_{\Phi_n} \sup_{f\in {\mathcal K}}\inf_{y\in B_{Y_n}}\|f-\Phi_n(y)\|_X,
\end{equation}
where the infimum is taken over all Lipschitz mappings
 $$
 \Phi_n:(B_{Y_n},\|\cdot\|_{Y_n})\to X,\quad B_{Y_n}:=\{y\in {\mathbb R}^n:\,\,\|y\|_{Y_n}\leq 1\},
 $$ 
that satisfy  
 the Lipschitz condition 
\begin{equation}
\nonumber
 \sup_{y,y'\in B_{Y_n}}
 \frac{\|\Phi_n(y)-\Phi_n(y')\|_X}{\|y-y'\|_{Y_n}}\leq \gamma,
\end{equation}
 with constant $\gamma$. 
We then define the
{\it Lipschitz} width
\begin{equation}
\nonumber
d_n^\gamma({\mathcal K})_X:=\inf_{k\leq n}\inf_{\|\cdot\|_{Y_k}}d^\gamma({\mathcal K},Y_k)_X,
\end{equation}
where the infimum is taken over  all norms $\|\cdot\|_{Y_k}$ in ${\mathbb R}^k$ and all $k\leq n$.    We observe the following analog to (\ref{hom})
\begin{equation}\label{homLip}
|t|d_n^{\gamma |t|}(t\mathcal{K})_X=d_n^{\gamma}(\mathcal{K})_X,  \mbox{ where } {t \mathcal K}:=\{tf:\,f\in \mathcal K\}.
\end{equation}

\section{Comparison between nonlinear Kolmogorov widths and  Lipschitz widths}
\label{S3}

In this section, 
we derive direct inequalities between the nonlinear Kolmogorov  widths  and the Lipschitz widths. We then use known relations between entropy numbers and Lipschitz widths to derive improvements of the (generalized) Carl's inequalities.

We first note the following comparison between the linear Kolmogorov widths and the Lipschitz widths, proven  in \cite{PW},
 see  Corollary 5.2.
\begin{theorem}
\label{T4}
For every $n\geq 1$ and  every compact set ${\mathcal K}\subset X$ we have
\begin{equation}
\nonumber
d^\gamma_n({\mathcal K})_X\leq d_n({\mathcal K})_X, \quad  \mbox{for every $\gamma\geq$}2\sup_{f\in {\mathcal K}}\|f\|_X.
\end{equation}
\end{theorem}

We next proceed with estimates between the nonlinear Kolmogorov width and the Lipschitz widths.
Clearly, it follows from the definition that 
$$ 
d_n({\mathcal K},N)_X\geq d_{nN}({\mathcal K})_X\geq d_{nN}^\gamma({\mathcal K})_X, \quad   \gamma= 2\sup_{f\in {\mathcal K}}\|f\|,
$$
where we have used in the last inequality the above theorem. Better estimates in the case of ${\mathcal K}$ being a subset 
of a Hilbert space $H$ or a general Banach space $X$ are described in the following lemmas.
 \begin{lemma}
\label{L1}
For every $n\geq 1$,  $N> 1$, and  every  compact  ${\mathcal K}$, subset of a Hilbert space $H$  such that $\sup_{f\in {\mathcal K}}\|f\|_H=1$, we have
\begin{equation}
\label{ll}
{  d^{(N+1)}_{n+1}({\mathcal K})_H}\leq d_n({\mathcal K},N)_H, \quad \hbox{and}\quad   d^{3}_{n+\lceil \log_2 N\rceil}({\mathcal K})_H\leq d_n({\mathcal K},N)_H.
\end{equation}
 \end{lemma}
\noindent
{\bf Proof:}  
Let us  fix $n,N\geq 1$, and  consider the $n$-dimensional linear spaces 
 $X_1,\dots,X_N$, $X_i\subset H$, $i=1,\ldots,N$. We define a norm $\|\cdot\|_{Y_{n+1}}$ on ${\mathbb R}^{n+1}$,
 $$
 \|(x,x_{n+1})\|_{Y_{n+1}}:=\max \left\{\|x\|_{\ell_2(\mathbb R^n)}, |x_{n+1}|\right\},\quad x:=(x_1,\ldots,x_n),
 $$
 whose unit ball is  
$$
{B}_{Y_{n+1}}:=\{(x, x_{n+1}) \  :   \ \|x\|_{\ell_2(\mathbb R^n)}\leq 1 \mbox{  and  } |x_{n+1}|\leq 1\}.
$$
Clearly
$$
B_{Y_{n+1}}=
B_{\ell_2(\mathbb R^n)}\times [-1,1], \quad \hbox{where}\quad B_{\ell_2(\mathbb R^n)}:= \{x\in \mathbb R^n \  :   \ \|x\|_{\ell_2(\mathbb R^n)}\leq 1\}.
$$
 We want to construct a Lipschitz mapping from $({B}_{Y_{n+1}}, \|\cdot\|_{Y_{n+1}})$ to $H$ whose image approximates well  $\mathcal K$.
We divide the interval $[-1,1]$ into $N$ subintervals $I_j$, $j=0,\ldots,N-1$, 
$$
 I_j:=[a_j,a_{j+1}], \quad a_j:=2j/N-1,
 $$
with  centers $c_j$ and consider the univariate continuous piecewise linear functions $\psi_j$, 
$\psi_j:([-1,1],|\cdot|)\rightarrow [0,1]$, $j=0,\ldots,N-1$, whose break points  are $\{a_0,\ldots,a_j,c_j,a_{j+1},\ldots,a_{N-1}\}$, and
$$
 \psi_j(c_j)=1, \quad \psi_j(a_k)=0, \quad k=0,\ldots,N-1.
$$
Let $(B_{X_j},\|\cdot\|_H)$ be the unit ball of the space $X_j\subset H$. We  fix an orthonormal  basis $\{\varphi^j_1,\ldots,\varphi^j_n\}$ in $X_j$ 
and consider the isometry map $\bar \psi_j$ from $B_{\ell_2(\mathbb R^n)}$ onto $B_{X_j}$,
$$
\bar \psi_j:(B_{\ell_2(\mathbb R^n)}, \|\cdot\|_{\ell_2({\mathbb R}^n)}) \to (B_{X_j},\|\cdot\|_H),
$$
defined as
\begin{equation}
\label{psi}
\bar \psi_j(x)=\bar \psi_j(x_1,\ldots,x_n):=\sum_{i=1}^nx_i\varphi^j_i.
\end{equation}
We use these  mappings to construct  $\Phi_{n+1}:({B}_{Y_{n+1}}, \|\cdot\|_{Y_{n+1}}) \rightarrow H$ as
$$
\Phi_{n+1}(x,x_{n+1}):=\sum_{j=0}^{N-1} \psi_j(x_{n+1}) \cdot \bar \psi_j(x).
$$
 Let us fix $(x,x_{n+1}),(x',x'_{n+1})\in B_{Y_{n+1}}$ and denote by 
$$
A:=\|\Phi_{n+1}(x,x_{n+1})-\Phi_{n+1}(x',x'_{n+1})\|_H.
$$
We want to derive an upper bound for $A$. Note that 
$\psi_j(x_{n+1})\neq 0$ if and only if $x_{n+1}\in I_j$.We 
 consider  the following two cases:

\begin{itemize}
\item  if  $x_{n+1},x'_{n+1}\in I_j$ for some $j=0,\ldots,N-1$, 
then  $\psi_j(x_{n+1})\neq 0$, $\psi_j(x'_{n+1})\neq 0$, $\psi_k(x_{n+1})=\psi_k(x'_{n+1})=0$ for all $k\neq j$, and therefore
\begin{eqnarray*}
A&=&\|\psi_j(x_{n+1})\bar \psi_j(x)-\psi_j(x'_{n+1})\bar \psi_j(x')\|_H\\
&\leq&
|\psi_j(x_{n+1})|\|\bar \psi_j(x)-\bar \psi_j(x')\|_H\\
\nonumber
&+&|\psi_j(x_{n+1}) -\psi_j(x'_{n+1})|\|\bar \psi_j(x')\|_H\\
\nonumber
&\leq& \|x-x'\|_{\ell_2({\mathbb R}^n)}+N|x_{n+1}-x'_{n+1}|\\
\label{P1}
&\leq &(N+1)\|(x,x_{n+1})-(x',x'_{n+1})\|_{Y_{n+1}}.
\end{eqnarray*}
\item if  $x_{n+1}\in I_j,x'_{n+1}\in I_k$ for some $j,k=0,\ldots,N-1$, $k\neq j$, we obtain that 
\begin{eqnarray*}
A=\|\psi_j(x_{n+1})\bar \psi_j(x)-\psi_k(x'_{n+1})\bar \psi_k(x')\|_H.
\end{eqnarray*}
We can assume without loss of generality that 
$$
x_{n+1}\leq a_{j+1}\leq a_k\leq x'_{n+1}.
$$
Since $\psi_j(a_{j+1})=\psi_k(a_k)=0$, we have
\begin{eqnarray*}
A&\leq&
\|\psi_j(x_{n+1})\bar \psi_j(x)-\psi_j(a_{j+1})\bar \psi_j(x)\|_H\\
&+&
\|\psi_k(a_k)\bar \psi_k(x)-\psi_k(x'_{n+1})\bar \psi_k(x')\|_H\\
&\leq &|\psi_j(x_{n+1}) -\psi_j(a_{j+1})|\|\bar \psi_j(x)\|_H\\
&+& \|\psi_k(a_k)\bar \psi_k(x)-\psi_k(x'_{n+1})\bar \psi_k(x')\|_H\\
&\leq &N|a_{j+1}-x_{n+1}|+ \|x-x'\|_{\ell_2({\mathbb R}^n)}+N|x'_{n+1}-a_k|\\
&\leq &N|x'_{n+1}-x_{n+1}|+ \|x-x'\|_{\ell_2({\mathbb R}^n)}\\
&\leq &(N+1)\|(x,x_{n+1})-(x',x'_{n+1})\|_{Y_{n+1}},
\end{eqnarray*}
where we have used arguments similar to the first case.
\end{itemize}

In both cases we have that 
$$
\|\Phi_{n+1}(x,x_{n+1})-\Phi_{n+1}(x',x'_{n+1})\|_H\leq (N+1)\|(x,x_{n+1})-(x',x'_{n+1})\|_{Y_{n+1}},
$$
and therefore $\Phi_{n+1}$ is an $(N+1)$-Lipschitz mapping.

 Since  
 $\sup_{f\in {\mathcal K}} \|f\|_H= 1$, the approximant $f_j$ to $f$ from $X_j$ will belong to 
$B_{X_j}$ since $f_j$ is the orthogonal projection of $f$ onto $X_j$.   Thus, it follows from the definition of $\bar \psi_j$ that 
there is $x^j\in B_{\ell_2(\mathbb R^n)}$, such that $\bar \psi_j(x^j)=f_j$, and therefore
$$
\Phi_{n+1}(x^j,c_j)=f_j, \quad \hbox{and}\quad \|f-f_j\|_H={\rm dist}(f,X_j)_H,
$$
which gives
 $$
 { d^{(N+1)}_{n+1}}({\mathcal K})_H\leq d_n({\mathcal K},N)_H.
 $$ 
To show the second part of (\ref{ll}), we determine $\ell\in {\mathbb N}$ such that 
$$
2^{\ell-1}<N\leq 2^\ell,
$$
and define a norm $\|\cdot\|_{{Y}_{n+\ell}}$ on ${\mathbb R}^{n+\ell}$ by
 $$
 \|(x,y)\|_{ Y_{n+\ell}}:=\max \left\{ \|x\|_{\ell_2(\mathbb R^n)}, \|y\|_{\ell_\infty({\mathbb R}^\ell)}\right\},
 $$
 where
 $$
 x:=(x_1,\ldots,x_n), \quad y:=(y_1,\ldots,y_\ell).
 $$
The  unit ball with respect to this norm is
$$
{B}_{Y_{n+\ell}}:=\{(x,y)\in \mathbb R^{n+\ell} \  :   \ \|x\|_{\ell_2(\mathbb R^n)}\leq 1 \mbox{  and  } \|y\|_{\ell_\infty({\mathbb R}^\ell)}\leq 1\}.
$$
  Like before, we have
$B_{Y_{n+\ell}}=B_{\ell_2(\mathbb R^n)}\times [-1,1]^\ell$.
Next, we consider the disjoint cubes $Q_j$, $j=1,\ldots,2^{\ell}$, of side length $1$ such that 
\begin{eqnarray*}
{\displaystyle [-1,1]^\ell=\cup_{j=1}^{2^\ell}Q_j}.
\end{eqnarray*}
We denote by
${\bf c}_j:=(c^j_{1},\dots,c^j_{\ell})\in\mathbb R^\ell$
the center of $Q_j$, $j=1,\ldots,2^\ell$,
and define the
functions $\phi_j:([-1,1]^\ell,\|\cdot \|_{\ell_\infty(\mathbb R^\ell)})\rightarrow [0,1]$ as 
$$
\phi_j(y):= 2\left ({\frac{1}{2}-\|{\bf c}_j-y\|_{\ell_\infty(\mathbb R^\ell)}}\right)_+, 
\quad j=1,\ldots,2^\ell,
$$
and $\Psi_{n+\ell}:({B}_{Y_{n+\ell}}, \|\cdot\|_{ Y_{n+\ell}}) \rightarrow H$ as
$$
\Psi_{n+\ell}(x,y):=\sum_{j=1}^{2^\ell} \phi_j(y) \cdot \bar \psi_j(x),
$$
where $\bar \psi_j$ are the mappings defined in (\ref{psi}).

Using the fact that for any two numbers $a,b$, we have $|a_+-b_+|\leq |a-b|$, we obtain that 
$$
|\phi_j(y)-\phi_j(y')|\leq 2|\|{\bf c}_j-y\|_{\ell_\infty(\mathbb R^\ell)}-\|{\bf c}_j-y'\|_{\ell_\infty(\mathbb R^\ell)}|\leq 2\|y-y'\|_{\ell_\infty(\mathbb R^\ell)}.
$$
Moreover, the supports of  the $\phi_j$'s are disjoint,  with $Q_j$ being the support of $\phi_j$, and $|\phi_j(y)|\leq 1$ for all $j$.
Now, following similar arguments as the ones for $\Phi_{n+1}$,  and denoting 
$$
B:=\|\Psi_{n+\ell}(x,y)-\Psi_{n+\ell}(x',y')\|_H,
$$
 we derive that:
\begin{itemize}
\item   if  $y,y'\in Q_j$ for some $j=1,\ldots,2^\ell$,
\begin{eqnarray*}
B=\|\phi_j(y)\bar\psi_j(x)-\phi_j(y')\bar\psi_j(x')\|_H
\leq3\|(x,y)-(x',y')\|_{Y_{n+\ell}}.
\end{eqnarray*}
\item if $y\in Q_j$ and $y'\in Q_k$, $k\neq j$, we consider the line segment 
$$
y+t(y'-y), \quad 0\leq t\leq 1,
$$ 
and   fix 
$$
d_j:=y+t_0(y'-y)\in \partial Q_j,  
$$
and
$$
b_k:=y+t_1(y'-y)\in \partial Q_k. 
$$
Clearly $t_0\leq t_1$, $\phi_j(d_j)=\phi_k(b_k)=0$, 
$$
\|y-d_j\|_{\ell_\infty(\mathbb R^\ell)}+
\|y'-b_k\|_{\ell_\infty(\mathbb R^\ell)}=(t_0+1-t_1)\|y-y'\|_{\ell_\infty(\mathbb R^\ell)}\leq
\|y-y'\|_{\ell_\infty(\mathbb R^\ell)},
$$
and similarly to the estimate for $A$, one obtains
\begin{eqnarray*}
B&=&\|\phi_j(y)\bar\psi_j(x)-\phi_k(y')\bar\psi_k(x')\|_H\\
&\leq&|\phi_j(y) -\phi_j(d_j)|\|\bar \psi_j(x)\|_H+\|\phi_k(b_k)\bar \psi_k(x)-\phi_k(y')\bar \psi_k(x')\|_H\\
&\leq&2\|d_j-y\|_{\ell_\infty(\mathbb R^\ell)}+ \|x-x'\|_{\ell_2({\mathbb R}^n)}+2\|y'-b_k\|_{\ell_\infty(\mathbb R^\ell)}\\
&\leq&2\|y-y'\|_{\ell_\infty(\mathbb R^\ell)}+ \|x-x'\|_{\ell_2({\mathbb R}^n)}\\
&\leq&3\|(x,y)-(x',y')\|_{Y_{n+\ell}}.
\end{eqnarray*}
\end{itemize}
Therefore, $\Psi_{n+\ell}$ is a
$3$-Lipschitz mapping. As before,  since
  $\sup_{f\in \mathcal K} \|f\|_H= 1$,  we obtain
 $$
 {d^{3}_{n+\lceil \log_2 N\rceil}(\mathcal K)_H}\leq d_n(\mathcal K,N)_H,
 $$
 where we have used  the fact that $\ell=\lceil \log_2 N\rceil$ and $\phi_j({\bf c}_j)=1$, $j=0, \ldots,N$.
The proof is completed. \hfill $\Box$
 
The case of arbitrary Banach space $X$ is based on the following lemma.

\begin{lemma}
\label{add}
Let $Y$ be an $n$-dimensional subspace of a Banach space $X$ and $(B_Y,\|\cdot\|_Y)$ be its unit ball.
Let $(B_Z,\|\cdot\|_H)$ be the unit ball in an $n$-dimensional subspace $Z$ of a Hilbert space $H$. Then,
 there exists a linear   map 
$$
\bar \psi:(B_Z,\|\cdot\|_H) \to Y,
$$
 with Lipschitz constant (i.e. norm ) at most $\sqrt n$ such that  $B_Y\subset \bar\psi (B_Z)$. In addition, if $X=L_p$, then the Lipschitz constant of $\bar \psi$ is at most $n^{|1/2-1/p|}$.
\end{lemma}
\noindent
{\bf Proof:} 
 It follows from the Fritz John theorem, see Chapter 3 in \cite{P} or \cite{KB}, that there exists an invertible linear operator 
 $\phi:(\mathbb R^n,\|\cdot\|_{\ell_2(\mathbb R^n)})\to Y$  onto $Y$ such that 
\begin{equation}
\phi(B_{\ell_2(\mathbb R^n)})\subset B_{ Y}\subset \sqrt{n}\phi(B_{\ell_2(\mathbb R^n)}).
\label{John} 
\end{equation}
Let us fix an orthonormal  basis $\varphi_1,\ldots,\varphi_n$ for $Z$ and consider the coordinate mapping
$\kappa_Z: Z \to \mathbb R^n   $ defined as
$$
\kappa_Z(g)=(x_1,\ldots,x_n)=x, \quad \hbox{where}\quad g=\sum_{j=1}^nx_j\varphi_j.
$$
This mapping  is isometry when $\mathbb R^n$ is equipped with the norm 
$$
\|x\|_{\ell_2(\mathbb R^n)}=\sqrt{\sum_{j=1}^n x_j^2} =\|g\|_Z.
$$
We now define the  linear mapping
$$
\tilde  \psi:=\phi \circ \kappa_Z:(Z,\|\cdot\|_H) \to Y,
$$
and notice that
$$
\tilde \psi(B_Z)\subset B_Y\subset \sqrt n  \tilde  \psi (B_Z).
$$
The first inclusion gives that $\tilde \psi$ has a norm (Lipschitz constant) $\leq 1$, and thus 
$\bar \psi :=\sqrt n \tilde \psi$ has a Lipschitz constant $\sqrt{n}$. 
The second inclusion shows that $B_Y\subset \bar\psi (B_Z)$, and therefore $\bar \psi$ is the desired mapping.
It follows from  \cite[Cor. 5]{ DRL}
that in the case of  $X=L_p$, we can replace $\sqrt n$  in (\ref{John}) by $n^{ |1/2-1/p|}$.
\hfill$\Box$

\begin{remark}
Note that since $\bar \psi$ is linear, we have that $\bar\psi(0)=0$, and  for every $z\in B_Z$,
\begin{equation}
\label{lin}
\|\bar \psi(z)\|_Y=\|\bar \psi(z)-\bar \psi(0)\|_Y\leq \sqrt{n}\|z\|_H\leq \sqrt{n},
\end{equation}
 where we can replace $\sqrt n$  by $n^{ |1/2-1/p|}$ in the case when $X=L_p$.
\end{remark}

 \begin{lemma}
\label{L2}
For every $n\geq 1$,  $N> 1$, and  every compact set ${\mathcal K}$ subset of a Banach space $X$
 with $\sup_{f\in {\mathcal K}}\|f\|_X=1$,  we have
\begin{equation}
\label{lltwo}
   d^{2(N+1)\sqrt{n}}_{n+1}({\mathcal K})_X\leq d_n({\mathcal K},N)_X, \quad \hbox{and}\quad 
  d^{ 6\sqrt{n}}_{n+\lceil \log_2 N\rceil}({\mathcal K})_X\leq d_n({\mathcal K},N)_X.
\end{equation}
 When $X=L_p$, we have
$$ 
 d^{2(N+1) n^{|1/2-1/p|}}_{n+1}({\mathcal K})_{L_p}\leq d_n({\mathcal K},N)_{L_p}, \quad \hbox{and}\quad   d^{ 6n^{|1/2-1/p|}}_{n+\lceil \log_2 N\rceil}({\mathcal K})_{L_p}\leq d_n({\mathcal K},N)_{L_p}.
$$
\end{lemma}
\noindent
{\bf Proof:}   
We  fix  $n$, $N> 1$, and  consider the $n$ dimensional linear spaces 
 $X_1,\dots,X_N$, $X_j\subset X$, $j=1,\ldots,N$, with  $(B_{X_j},\|\cdot\|_X)$ being the unit ball of $X_j$. 
 For  a fixed  $j=1,\ldots,N$, we apply Lemma \ref{add} with $Y=X_j$ and $Z=\ell_2(\mathbb R^n)$  to find an 
  $M$-Lipschitz mapping $ \bar\Psi_j$, where $M=\sqrt n$ or $n^{|1/p-1/2|}$, depending on whether $X$ is a general Banach space or $L_p$, such that
\begin{equation}
\label{lk}
 \bar\Psi_j:(B_{\ell_2(\mathbb R^n)}, \|\cdot\|_{\ell_2({\mathbb R}^n)}) \to X_j, \quad \hbox{and}\quad
 B_{X_j}\subset \bar\Psi_j(B_{\ell_2(\mathbb R^n)}).
\end{equation}

We show (\ref{lltwo})  by proceeding  as in the proof of Lemma \ref{L1} 
 and defining a mapping 
   $\Theta_{n+1}:({B}_{Y_{n+1}}, \|\cdot\|_{Y_{n+1}}) \rightarrow X$ as
$$
\Theta_{n+1}(x,x_{n+1}):= 2\sum_{j=0}^{N-1} \psi_j(x_{n+1}) \cdot \bar \Psi_j(x),
$$
where $\psi_j$ and $({B}_{Y_{n+1}}, \|\cdot\|_{Y_{n+1}})$ are as in Lemma \ref{L1}. We fix
$(x,x_{n+1})$,  $(x',x'_{n+1})$, denote by 
$$
C:=\|\Theta_{n+1}(x,x_{n+1})-\Theta_{n+1}(x',x'_{n+1})\|_X,
$$
and show in a similar way that 
\begin{itemize}
\item if $x_{n+1},x'_{n+1}\in I_j$ for some $j=0,\ldots,N-1$,
\begin{eqnarray*}
\frac{C}{2}&=&\|\psi_j(x_{n+1})\bar \Psi_j(x)-\psi_j(x'_{n+1})\bar \Psi_j(x')\|_X\\
&\leq&|\psi_j(x_{n+1})|\|\bar \Psi_j(x)-\bar \Psi_j(x')\|_X\\
&+&|\psi_j(x_{n+1}) -\psi_j(x'_{n+1})|\|\bar \Psi_j(x')\|_X\\
&\leq&  M\|x-x'\|_{\ell_2({\mathbb R}^n)}+NM|x_{n+1}-x'_{n+1}|\\
&\leq & M (N+1)\|(x,x_{n+1})-(x',x'_{n+1})\|_{Y_{n+1}},
\end{eqnarray*}
 where we have used (\ref{lin}).
\item if  $x_{n+1}\in I_j,x'_{n+1}\in I_k$ for some $j,k=0,\ldots,N-1$, $k\neq j$,
\begin{eqnarray*}
\frac{C}{2}&\leq&
\|\psi_j(x_{n+1})\bar \Psi_j(x)-\psi_j(a_{j+1})\bar \Psi_j(x)\|_X\\
&+&
\|\psi_k(a_k)\bar \Psi_k(x)-\psi_k(x'_{n+1})\bar \Psi_k(x')\|_X\\
&\leq &|\psi_j(x_{n+1}) -\psi_j(a_{j+1})|\|\bar \Psi_j(x)\|_X\\
&+& \|\psi_k(a_k)\bar \Psi_k(x)-\psi_k(x'_{n+1})\bar \Psi_k(x')\|_X\\
&\leq &NM|a_{j+1}-x_{n+1}|{ +  M}\|x-x'\|_{\ell_2({\mathbb R}^n)}+NM|x'_{n+1}-a_k|\\
&\leq & NM|x'_{n+1}-x_{n+1}|{ + M}\|x-x'\|_{\ell_2({\mathbb R}^n)}\\
&\leq &{ M(N+1)}\|(x,x_{n+1})-(x',x'_{n+1})\|_{Y_{n+1}}.
\end{eqnarray*}
\end{itemize}
In conclusion,
$$
\|\Theta_{n+1}(x,x_{n+1})-\Theta_{n+1}(x',x'_{n+1})\|_H\leq { 2M(N+1)}\|(x,x_{n+1})-(x',x'_{n+1})\|_{Y_{n+1}},
$$
and therefore $\Theta_{n+1}$ is a { $2M(N+1)$-Lipschitz mapping.}

Note that if $f_j$ is the approximant  to $f$ from $X_j$, then
\begin{equation}\label{approximant}
\|f-f_j\|_X\leq \|f\|_X\quad \Rightarrow \quad \|f_j\|_X\leq \|f-f_j\|_X+\|f\|_X\leq 2\|f\|_X\leq 2,
\end{equation}
where we have used that   $\sup_{f\in {\mathcal K}} \|f\|_X= 1$. Thus $f_j\in 2 B_{X_j}$.
   It follows from Lemma \ref{add} that since $B_{X_j}\subset \bar \Psi_j(B_{\ell_2(\mathbb R^n)})$, 
 there is $x^j\in B_{\ell_2(\mathbb R^n)}$, such that $\bar \Psi_j(x^j)=\frac{1}{2}f_j$. Therefore
$$
\Theta_{n+1}(x^j,c_j)=f_j, \quad \hbox{and}\quad \|f-f_j\|_X={\rm dist}(f,X_j)_X,
$$
which gives
 $$
 { d^{2M(N+1)}_{n+1}}({\mathcal K})_{X}\leq d_n({\mathcal K},N)_X.
 $$
  
To show the second part of (\ref{lltwo}), we define
 $\Xi_{n+\ell}:({B}_{ Y_{n+\ell}}, \|\cdot\|_{ Y_{n+\ell}}) \rightarrow X$ as
$$
\Xi_{n+\ell}(x,y):=2\sum_{j=1}^{2^\ell} \phi_j(y) \cdot \bar \Psi_j(x),
$$
where $\phi_j$ and $({B}_{Y_{n+\ell}}, \|\cdot\|_{ Y_{n+\ell}}) $ are the same as in Lemma \ref{L1}
and $\bar \Psi_j$ is defined in (\ref{lk}).
For fixed $(x,y), (x',y')\in {B}_{ Y_{n+\ell}}$, we denote by
$$
D:=\|\Xi_{n+\ell}(x,y)-\Xi_{n+\ell}(x',y')\|_X
$$
and consider the following cases
\begin{itemize}
\item   if  $y,y'\in Q_j$ for some $j=1,\ldots,2^\ell$, we have
\begin{eqnarray*}
\frac{D}{2}\leq  3M\|(x,y)-(x',y')\|_{ Y_{n+\ell}}.
\end{eqnarray*}
\item if $y\in Q_j$ and $y'\in Q_k$, $k\neq j$,  similarly to the estimate for C, we obtain
\begin{eqnarray*}
{\frac{D}{2}}&=&\|\phi_j(y)\bar\Psi_j(x)-\phi_k(y')\bar\Psi_k(x')\|_X\\
&\leq&|\phi_j(y) -\phi_j(d_j)|\|\bar \Psi_j(x)\|_X+\|\phi_k(b_k)\bar \psi_k(x)-\phi_k(y')\bar \Psi_k(x')\|_X\\
&\leq&2M\|d_j-y\|_{\ell_\infty(\mathbb R^\ell)}+ { M}\|x-x'\|_{\ell_2({\mathbb R}^n)}+2M\|y'-b_k\|_{\ell_\infty(\mathbb R^\ell)}\\
&\leq&2M\|y-y'\|_{\ell_\infty(\mathbb R^\ell)}+{ M} \|x-x'\|_{\ell_2({\mathbb R}^n)}\\
&\leq&{ 3M}\|(x,y)-(x',y')\|_{ Y_{n+\ell}}.
\end{eqnarray*}
\end{itemize}
The latter estimate implies that $\Xi_{n+\ell}$ is a
{ $6M$}-Lipschitz mapping, and since   $\sup_{f\in \mathcal K} \|f\|_X= 1$,  we obtain
 $$
 d^{6M}_{n+\lceil \log_2 N\rceil}(\mathcal K)_X\leq d_n(\mathcal K,N)_X.
 $$
The proof is completed. \hfill $\Box$

{ \begin{remark} Note  that Lemma \ref{L2} with $X=L_2$ can be used instead of Lemma \ref{L1}. However, we 
have decided to present both lemmas since  better Lipschitz constants are obtained when working directly with a Hilbert space $H$.

\end{remark}
}

\begin{remark}
It follows from (\ref{hom}) and (\ref{homLip})  that lemmas similar to Lemma \ref{L1} and Lemma \ref{L2}
can be stated in the case when $\sup_{f\in {\mathcal K}}\|f\|_H\neq 1$, or $\sup_{f\in {\mathcal K}}\|f\|_X\neq 1$, respectively.
 \end{remark}

\section{Main results}
\label{S4}

In this section, we provide estimates from above and below that connect the behavior of the linear and nonlinear Kolmogorov widths of $\mathcal K$ with its entropy numbers.
In what follows we assume that
$\sup_{f\in \mathcal K} \|f\|_H=1$ in the case of Hilbert space, or  $\sup_{f\in \mathcal K} \|f\|_X=1$ in the case of a general Banach space. Similar results hold  if this supremum is not 1.

Our approach  of deriving estimates from below utilizes  some known results  for Lipschitz widths stated below, see Theorem 4.7 in \cite{PW}.
\begin{theorem}
\label{T1}
Let  ${\mathcal K}\subset X$ be a  compact subset of a Banach space $X$,  $n\in N$, and $d_n^{\gamma}({\mathcal K})_X$ be the   Lipschitz width for ${\mathcal K}$ with Lipschitz constant  
$\gamma\geq 2{\rm rad}({\mathcal K})$. Then the following holds:
\begin{enumerate}
\item If for   $\alpha>0$, $\beta\in R$ and a constant $C>0$, we have 
\begin{eqnarray}
\nonumber
e_n({\mathcal K})_X\geq C \frac{[\log_2n]^\beta}{n^\alpha},\quad n=1,2,\ldots,\quad \quad \hbox{then}\quad 
d_n^{\gamma}({\mathcal K})_X\geq C'\frac{[\log_2n]^\beta}{n^{\alpha}[\log_2 n]^{\alpha}},
\end{eqnarray}
for $n=1,2,\ldots$, where  $C'>0$ is a fixed constant.
\item If for   $\alpha>0$ and $C>0$,  we have 
\begin{equation}
\nonumber
e_n({\mathcal K})_X\geq { C}\frac{1}{[\log_2 n]^\alpha},\quad n=1,2,\ldots,\quad \quad \hbox{then}\quad 
d_n^{\gamma}({\mathcal K})_X\geq \frac{C'}{[\log_2 n]^{\alpha}},
\end{equation}
for $n=1,2,\ldots$, where  $C'>0$ is a fixed constant.
\item If for   $0<\alpha<1$ and $C,c>0$,  we have 
\begin{eqnarray}
\nonumber
e_n({\mathcal K})_X\geq C 2^{-cn^{\alpha}},\quad n=1,2,\ldots,\quad \quad \hbox{then}\quad 
d_n^{\gamma}({\mathcal K})_X\geq C'2^{-c'n^{\alpha/(1-\alpha)}},
\end{eqnarray}
for $n=1,2,\ldots$, where  $C',c'>0$, are fixed constants.
\end{enumerate}
\end{theorem}

\subsection{Estimates from below for the linear Kolmogorov width.}
The above theorem, combined with  Theorem \ref{T4},  gives the following relations between linear Kolmogorov widths and 
entropy numbers.
\begin{theorem}
\label{T1new}
Let  ${\mathcal K}\subset X$ be a  compact subset of a Banach space $X$,  $n\in \mathbb N$, and $d_n({\mathcal K})_X$ be the   $n$-th linear Kolmogorov width for ${\mathcal K}$. Then the following holds:
\begin{enumerate}
\item If for   $\alpha>0$, $\beta\in \mathbb{ R}$, $C>0$, we have 
\begin{eqnarray}
\nonumber
e_n({\mathcal K})_X\geq C \frac{[\log_2n]^\beta}{n^\alpha},\quad n=1,2,\ldots,\quad \quad \hbox{then}\quad 
d_n({\mathcal K})_X\geq C'\frac{[\log_2n]^\beta}{n^{\alpha}[\log_2 n]^{\alpha}},
\end{eqnarray}
for $n=1,2,\ldots$, where  $C'>0$ is a fixed constant.
\item If for   $\alpha>0$,  $C>0$, we have 
\begin{equation}
\nonumber
e_n({\mathcal K})_X{ \geq \frac{C}{[\log_2 n]^\alpha}} ,\quad n=1,2,\ldots,\quad \quad \hbox{then}\quad 
d_n({\mathcal K})_X\geq C'\frac{1}{[\log_2 n]^\alpha},
\end{equation}
for $n=1,2,\ldots$, where  $C'>0$ is a fixed constant.
\item If for   $0<\alpha<1$,  $C,c>0$ we have 
\begin{eqnarray}
\nonumber
e_n({\mathcal K})_X\geq C 2^{-cn^{\alpha}},\quad n=1,2,\ldots,\quad \quad\hbox{then}\quad 
d_n({\mathcal K})_X\geq C'2^{-c'n^{\alpha/(1-\alpha)}},
\end{eqnarray}
for $n=1,2,\ldots$, where  $C',c'>0$ are fixed constants.
\end{enumerate}
\end{theorem}
\noindent
{\bf Proof:} The statement follows from Theorem \ref{T4}, Theorem \ref{T1} and the inequality $\sup_{f\in K}\|f\|_X\geq {\rm rad}({\mathcal K})$. 
\hfill $\Box$

\subsection{Estimates from below for the nonlinear Kolmogorov width, the Hilbert space case.}
Using Lemma \ref{L1} and Theorem \ref{T1}, we obtain similar estimates for  $d_{n-1}({\mathcal K},N)_H$.
\begin{theorem}
\label{T2new}
Let  ${\mathcal K}\subset H$ be a  compact subset of a Hilbert space $H$ and $d_n(\mathcal K, N)_H$, $n\in \mathbb N$, $N>1$,  be the   nonlinear Kolmogorov width for ${\mathcal K}$. Then the following holds:
\begin{itemize}
\item If for    $\alpha>0$, $\beta\in \mathbb R$, and  $C>0$  the entropy numbers satisfy $e_n({\mathcal K})_H\geq C \frac{[\log_2n]^\beta}{n^\alpha}$, $n=1,2,\ldots$, then  there is a constant $C''>0$ such that for every $N> 1$we have
\begin{eqnarray}
\label{main12}
d_{n-1}({\mathcal K},N)_H\geq C''\frac{[\log_2(n+\lceil\log_2 N\rceil)]^{\beta-\alpha}}{[n+\lceil\log_2 N\rceil]^{\alpha}},\quad n=1,2,\ldots.
\end{eqnarray}
\item If for    $\alpha>0$ and $C>0$,  the entropy numbers satisfy the inequality $e_n({\mathcal K})_H\geq \frac{C}{[\log_2 n]^\alpha}$, $n=1,2,\ldots$, then there is a constant $C''>0$ such that for every $N> 1$ we have
\begin{equation}
\label{main22}
d_{n-1}({\mathcal K},N)_H\geq C'' \frac{1}{[\log_2 (n+\lceil\log_2 N\rceil)]^\alpha},\quad n=1,2,\ldots.
\end{equation}

\item If for   $0<\alpha<1$ and $C,c>0$,   the entropy numbers satisfy the inequality $e_n({\mathcal K})_H\geq C 2^{-cn^{\alpha}}$,  $n=1,2,\ldots$, then there are constants $C'',c''>0$ such that  for every $N>1$
\begin{eqnarray}
\label{main33}
d_{n-1}({\mathcal K},N)_H\geq C''2^{-c''(n+\lceil\log_2 N\rceil)^{\alpha/(1-\alpha)}},\quad n=1,2,\ldots.
\end{eqnarray}
\end{itemize}
\end{theorem}
\noindent
{\bf Proof:} To show  (\ref{main12}), we apply Lemma \ref{L1}, Theorem \ref{T1} with  a value $\gamma=\max\{2{\rm rad}(\mathcal K), 3\}$,  and use the monotonicity of the Lipschitz width as a function of $\gamma$ to derive that 
\begin{eqnarray*}
d_{n-1}({\mathcal K},N)_H&\geq& { d_{n+\lceil \log_2 N\rceil}^{3}(\mathcal K)_H}\\
&\geq &d_{n+\lceil \log_2 N\rceil}^{\gamma}(\mathcal K)_H\geq
C \frac{[\log_2(n+\lceil \log_2 N\rceil)]^{\beta-\alpha}}{[{n+\lceil \log_2 N\rceil}]^\alpha}.
\end{eqnarray*}
We omit the  proof of the rest of the theorem since it is similar to the case  already discussed. 
\hfill$\Box$

Note that the above theorem holds for any value of $N$. In the cases when  $N=\lambda^n$, with $\lambda>1$,
or $N=n^{an}$, with $a>0$, we obtain two corollaries.
\begin{corollary}
\label{c1}
Let $\mathcal K\subset H$ be a compact subset of a Hilbert space $H$. Then the following holds:
\begin{itemize}
\item If $e_n({\mathcal K})_H\geq C \frac{[\log_2n]^\beta}{n^\alpha},\quad n=1,2,\ldots$, then
$$
d_{n-1}({\mathcal K},\lambda^n)_H\geq C''\frac{[\log_2 n]^{\beta-\alpha}}{n^{\alpha}},\quad n=2,3,\ldots.
$$
\item If $e_n({\mathcal K})_H\geq C\frac{1}{[\log_2 n]^\alpha},\quad n=1,2,\ldots$, then
$$
d_{n-1}({\mathcal K},\lambda^n)_H\geq C''\frac{1}{[\log_2 n]^\alpha},\quad n=2,3,\ldots.
$$
\item If $e_n({\mathcal K})_H\geq C 2^{-cn^{\alpha}},\quad n=1,2,\ldots$, then
$$
d_{n-1}({\mathcal K},\lambda^n)_H\geq C''2^{-c''n^{\alpha/(1-\alpha)}},\quad n=2,3,\ldots.
$$
\end{itemize}
\end{corollary}

\begin{corollary}
\label{c2}
Let $\mathcal K\subset H$ be a compact subset of a Hilbert space $H$. Then the following holds:
\begin{itemize}
\item If $e_n({\mathcal K})_H\geq C \frac{[\log_2n]^\beta}{n^\alpha},\quad n=1,2,\ldots$ then
$$
d_{n-1}({\mathcal K},n^{an})_H\geq C''\frac{[\log_2 n]^{\beta-2\alpha}}{n^{\alpha}},\quad n=2,3,\ldots.
$$
\item If $e_n({\mathcal K})_H\geq C\frac{1}{[\log_2 n]^\alpha},\quad n=1,2,\ldots$, then
$$
d_{n-1}({\mathcal K},n^{an})_H\geq C''\frac{1}{[\log_2 n]^\alpha},\quad n=2,3,\ldots.
$$
\item If $e_n({\mathcal K})_H\geq C 2^{-cn^{\alpha}},\quad n=1,2,\ldots$, then
$$
d_{n-1}({\mathcal K},n^{an})_H\geq C''2^{-c''[n\log_2 n]^{\alpha/(1-\alpha)}},\quad n=2,3,\ldots.
$$
\end{itemize}
\end{corollary}
{
\noindent
{\bf Proof:} 
We outline the proof of only the first statement. It follows from (\ref{main12}) with $N=n^{an}$ that
\begin{eqnarray*}
d_{n-1}({\mathcal K},n^{an})_H&\geq& C''\frac{[\log_2(n+an\log_2 n)]^{\beta-\alpha}}{[n+an\log_2 n]^{\alpha}}
\geq C_1 \frac{[\log_2(n+an\log_2 n)]^{\beta-\alpha}}{[n\log_2 n]^{\alpha}}\\
&\geq&C_2\frac{[\log_2 n]^{\beta-\alpha}}{[n\log_2 n]^{\alpha}},
\end{eqnarray*}
where we have used that for $n$ big enough
$$
\log_2 n\leq \log_2(n+an\log_2 n)\leq 2\log_2 n.
$$
\hfill $\Box$
\subsubsection{Examples} 
Here,  we provide an example  which shows that some of  the estimates in Corollary \ref{c1} are sharp.
We consider the Hilbert space $\ell_2:=\{x=(x_1,x_2,\ldots):\,\,\sum_{j=1}^\infty|x_j|^2<\infty\}$ with a standard basis 
$\{e_j)_{j=1}^\infty$ and the strictly decreasing sequence  $\sigma=\{\sigma_j\}_{j=1}^\infty$ of positive numbers $\sigma_j$ which converge to $0$ with $\sigma_1=1$. We then define the compact set 
 $$
 {\mathcal K}_\sigma:=\{\sigma_j e_j\}_{j=1}^\infty\cup\{0\}\subset \ell_2
 $$
and prove the following lemma.
\begin{lemma}
\label{Lexp}
Every set $\mathcal K_\sigma\subset \ell_2$  has inner entropy numbers
$$
\tilde e_n({\mathcal K}_\sigma)_{\ell_2}=\sqrt{\sigma_{2^n}^2+\sigma_{2^n+1}^2}, \quad n=1,2,\ldots,
$$ 
and nonlinear Kolmogorov width
$$
d_n({\mathcal K}_\sigma, N)_{\ell_2}\leq \sigma_{nN+1}, \quad N>1, \quad n=1,2,\ldots.
$$
\end{lemma}
\noindent
{\bf Proof:} 
Since
$$
 \|\sigma_j e_j -\sigma_{j^\prime}e_{j^\prime}\|_{\ell_2}= \sqrt{ \sigma_j^2+\sigma_{j^\prime}^2}\leq \sqrt {\sigma_{j}^2+\sigma_{j+1}^2}, \quad \hbox{for all}\quad j'\geq j+1,
$$
and
$$
\|\sigma_j e_j -0\|_{\ell_2}= \sigma_j<\sqrt{\sigma_{j}^2+\sigma_{j+1}^2},
$$
we have that the ball with center $\sigma_j e_j$ and radius  $r_j:=\sqrt{\sigma_j^2+\sigma_{j+1}^2}$   contains $0$ and all points $\sigma_{j^\prime} e_{j^\prime}$ with $j^\prime>j$, but  none of the points $\sigma_{j^\prime} e_{j^\prime}$ with $j^\prime <j$. Thus, if we look for $2^n$ balls with centers in $\mathcal K_\sigma$,  covering $\mathcal K_\sigma$, and with smallest radius,  these are the balls $B(\sigma_j e_j, r_{2^n})$,  $j=1,2,\dots, 2^n$, with centers $\sigma_je_j$ and radius $r_{2^n}$. The $j$-th  ball does not contain the first $(j-1)$ points $\sigma_{j'}e_{j'}$, $1\leq j'\leq j-1$, from 
$\mathcal K_\sigma$, but contains the rest of the points $\{\sigma_{i} e_i\}_{i=j}^\infty \cup \{0\}$. Therefore, we have that 
$$
\tilde e_n(\mathcal K_\sigma)_{\ell_2}= r_{2^n}.
$$
 To prove the second statement,  we define the $n$-dimensional spaces 
 $$
 X_s:={\rm span}\{e_j\}_{(s-1)n+1}^{sn}, \quad s=1,2,\dots,N.
 $$
Clearly $0,\sigma_je_j\in \bigcup_{s=1}^N X_s$  for $j=1,\dots,nN$, and for $j>nN$ we have 
$$
{\rm dist}(\sigma_j e_j,\bigcup_{s=1}^N X_s)_{\ell_2}
=\sigma_j.
$$
Thus, $d_n({\mathcal K}_\sigma, N)_{\ell_2}\leq \sigma_{nN+1}$, and  the proof is completed. 
\hfill $\Box$

For our particular example we fix $\alpha>0$, select the sequence $\{\sigma_j\}_{j=1}^\infty$ to be 
\begin{equation}
\label{exp}
\sigma_j=\frac{1}{[\log_2 \log_2(j+3)]^\alpha},\quad j=1,2,\dots,
\end{equation}
and show in the following lemma that the estimate in Corollary \ref{c1} cannot be improved.

\begin{lemma}
\label{Lexp1}
The set $\mathcal K:=\mathcal K_\sigma$ defined by the sequence {\rm (\ref{exp})} has the following properties:
$$
e_n(\mathcal K)_{\ell_2} \asymp  (\log_2 n )^{-\alpha}, \quad \hbox{and} \quad d_{n-1}(\mathcal K, \lambda^n)_{\ell_2}
\asymp  (\log_2 n )^{-\alpha}, \quad n=2,3, \ldots.
$$
\end{lemma}
\noindent
{\bf Proof:} 
 It follows from (\ref{innerentropy}) and Lemma 
\ref{Lexp} that 
$$ 
e_n(\mathcal K)_{\ell_2} \asymp \sigma_{2^n}\asymp (\log_2 n )^{-\alpha},
$$
and that
\begin{eqnarray*}
 d_{n-1}(\mathcal K, \lambda^n)_{\ell_2}&\leq& \sigma_{(n-1)\lambda^n +1}=\frac{1}{[\log_2 \log_2((n-1)
\lambda^n+1)]^\alpha}  \\
&\leq &\frac{C}{(\log_2 \log_2 \lambda^n)^{\alpha}}=\frac{C}{(  \log_2 n +\log_2 \log_2 \lambda)^\alpha}\\
&\leq&\frac{C'}{(\log_2 n)^\alpha}.
\end{eqnarray*}
The estimate from below follows from Corollary \ref{c1}.
\hfill $\Box$
}

\subsection{Estimates from below for the nonlinear Kolmogorov width, the Banach space case.}

To prove an estimate from below in the Banach space case, we use the following statement from \cite{PW}, see Theorem 7.3 in \cite{PW}.

\begin{theorem}
\label{TT}
Let   $\mathcal K\subset X$  be a compact subset of a Banach space $X$. Consider the Lipschitz width $d_n^{\gamma_n}(\mathcal K)_X$ with 
$\gamma_n=cn^{\delta}\lambda^n$, $\delta\in\mathbb R$,  $\lambda >1$, and $c>0$. If for some constants $c_1>0,\alpha>0$  we have  
$
e_n(\mathcal K) _X> c_1 (\log_2 n)^{-\alpha},$  $n=1,2,\dots,
$
then there exists a constant $C>0$ such that 
\begin{equation} 
\nonumber
d_n^{\gamma_n}(\mathcal K)_X\geq C(\log_2 n)^{-\alpha}, \quad n=1,2,\dots .
\end{equation}
\end{theorem}

We now use Lemma \ref{L2} and the above statement to prove the following theorem.
\begin{theorem}
\label{T3new}
Let  ${\mathcal K}\subset X$ be a  compact subset of a Banach space $X$ and $d_n(\mathcal K,N)_X$, $n\in \mathbb N$, $N>1$,  be the   nonlinear Kolmogorov width for ${\mathcal K}$.  If there is  $\alpha>0$ and $C>0$ such that  the entropy numbers 
$e_n({\mathcal K})_X\geq C \frac{1}{[\log_2 n]^\alpha}$, $n=1,2,\ldots,$
 then  there is an absolute constant $C''>0$ such that 
\begin{equation}
\nonumber
d_{n-1}({\mathcal K},\lambda^n)_X\geq C'' \frac{1}{[\log_2 n]^\alpha},\quad n=2,\ldots .
\end{equation}
\end{theorem}
\noindent
{\bf Proof:} 
We apply Lemma \ref{L2}, Theorem \ref{TT} with $\gamma=2(\lambda^n+1)\sqrt{n}$  and use the monotonicity of the Lipschitz width as a function of $\gamma$ to derive that 
$$
d_{n-1}({\mathcal K},\lambda^n)_X\geq d_n^{2(\lambda^n+1)\sqrt{n}}(\mathcal K)_X\geq d_n^{c\sqrt{n}\lambda^n}(\mathcal K)_X\geq
C'\frac{1}{[\log_2 n]^\alpha}.
$$
\hfill$\Box$

\subsection{Estimates from above for the entropy numbers}
The next proposition provides us with a tool to derive estimates for the entropy numbers of $\mathcal K$ if we have 
a knowledge about the behavior of the nonlinear Kolmogorov widths $d_n(\mathcal K,N)_X$.

\begin{proposition}
\label{Tem1}
Let  ${\mathcal K}\subset X$ with ${\rm rad}(\mathcal{K})<1$ be a  compact subset of a Banach space $X$ and $d_n(\mathcal K,N)_X$, $N>1$,  $n\in \mathbb N$, be the   nonlinear Kolmogorov width for ${\mathcal K}$. If  for some $1>\epsilon>0$
we have  $d_n(\mathcal{K},N)_X<\epsilon$,  then there exists an absolute constant $c>0$ such that $\tilde P_{3\epsilon}(\mathcal{K}) \leq N (c/\epsilon)^n$ and
$$
 e_{\lceil \log_2 \mu \rceil}(\mathcal{K})_X\leq 3cN^{1/n}\mu^{-1/n},\quad  \hbox{with}\quad \mu=\tilde P_{3\epsilon}(\mathcal{K}).
$$
\end{proposition}
\noindent
{\bf Proof:} 
Since $d_n(\mathcal{K},N)_X<\epsilon$, it follows from the definition of the nonlinear Kolmogorov width that 
there exist $n$-dimensional subspaces $X_j\subset X$, $j=1,\dots,N$, for which
$$
\sup_{f\in \mathcal K}\,\,\,\inf_{X_j, \,j=1,\ldots,N}{\rm dist}(f,X_j)_X<\epsilon.
$$
Let $\{k_j\}_{j=1}^\mu\subset \mathcal{K}$ be a maximal $3\epsilon$-packing in $\mathcal{K}$, i.e. $\mu=\tilde P_{3\epsilon}(\mathcal{K}) $.
  Then, for each $k_i$ there exists $x_i\in X_{j(i)}\subset \bigcup_{j=1}^N X_j$ such that $\|k_i-x_i\|_X<\epsilon$, $i=1,\dots,\mu$, and  we can estimate the difference $\|x_i-x_{i^\prime}\|_X$, $i\neq i'$,
\begin{equation}
\label{qq}
 \|x_i-x_{i^\prime}\|_X\geq \|k_i-k_{i^\prime}\|_X-\|x_{i^\prime}-k_{i^\prime}\|_X-\|x_i-k_{i}\|_X> \epsilon.
\end{equation}
The condition  ${\rm rad}(\mathcal{K})<1$ implies that there exists $y\in X$ such that 
$$
\mathcal{K}\subset B(y,1)_X,
$$
and therefore
$$
\|x_i-y\|_X \leq \|x_i-k_i\|_X+\|k_i-y\|_X < 1+\epsilon.
$$
 Let $y_j$ be the closest to $y$  element from $X_{j(i)}\subset \bigcup_{j=1}^N X_j$. Then for any 
 $x_i\in \{x_1,\ldots,x_N\}$ ($x_1,\ldots,x_N$ may not be necessary all different)
 we have 
 $$
 \|y_j-x_i\|_X\leq \|y_j-y\|_X+\|y-x_i\|_X\leq 2\|y-x_i\|_X< 2(1+\epsilon),
 $$
which leads to
$$ 
\{x_i\}_{i=1}^\mu \subset \bigcup_{j=1}^N B_{X_{j(i)}}(y_j,2(1+\epsilon))\subset \bigcup_{j=1}^N B_{X_{j(i)}}(y_j,4)
, \quad \hbox{where}\quad 
y_j\in X_{j(i)}.
$$
It follows from (\ref{qq}) that the set $\{x_i\}_{i=1}^\mu$ is an $\epsilon$-packing for $\bigcup_{j=1}^N B_{X_{j(i)}}(y_j,4)$. 
We next use (\ref{known})   to derive
$$
\mu\leq \widetilde{ \mathcal N}_{\epsilon/2}\left(\bigcup_{j=1}^N B_{X_{j(i)}}(y_j,4)\right)\leq N(c/\epsilon)^n,
$$
where $c>0$ is an absolute constant. Note that  we have applied  the inequality 
$$
\widetilde{ \mathcal N}_{\epsilon/2}(B_{X_{j(i)}}(y_j,4))\leq (c/\epsilon)^n,\quad j=1,\ldots,N,
$$
from \cite[Chp. 15 Prop.1.3]{LGM}. 
In terms of entropy numbers we can write 
$$
e_{\lceil \log_2 \mu \rceil}(\mathcal{K})_X\leq 3 \epsilon\leq 3cN^{1/n}\mu^{-1/n}.
$$
\hfill $\Box$

We use Proposition \ref{Tem1} to obtain estimates from above for the entropy numbers $e_m(\mathcal K)_X$ of $\mathcal K$.
{ A similar estimate but for a different range of $m$ and  some specific  values  of $N$ has been recently presented in \cite{T1}.}
\begin{lemma}
\label{L6}
Let $\mathcal K\subset X$ be a compact subset of a Banach space $X$ with ${\rm rad}(\mathcal{K})<1$.
If for $\alpha>0$, $\beta\in \mathbb R$, $\lambda>1$, and $c_0>0$ we have  that 
$$
d_n(\mathcal{K},\lambda^n)_X\leq c_0\frac{[\log_2 n]^\beta}{ n^{\alpha}},
$$
for some $n>n_0(c_0,\alpha,\beta,\lambda)$,
then 
$$
e_m(\mathcal K)_X<C\frac{[\log_2 m]^{\alpha+\beta}}{ m^{\alpha}}, \quad \hbox{with}\quad m=2\alpha n\log_2n,
$$
where 
$C$ is a  fixed constant depending only on $\lambda,\alpha,\beta,c_0$.
  \end{lemma}
\noindent
{\bf Proof:}   It follows from Proposition \ref{Tem1} with $\epsilon=c_0\frac{ [\log_2 n]^\beta}{ n^{\alpha}}$ 
that 
$$
\log_2 \mu\leq n[\log_2(\lambda c)+\alpha\log_2n-\log_2 c_0-\beta\log_2(\log_2 n)]\leq 2\alpha n\log_2 n,
$$
 for $n>n_0$ where $n_0$ depends only $c_0$, $\lambda$, $\alpha$ and $\beta$. For such $n$'s we have
 $$
e_{2\alpha n \log_2 n}(\mathcal{K})_X \leq \frac{3c_0[\log_2 n]^\beta}{n^\alpha}.
 $$
Setting $m=2\alpha n \log_2 n$ gives
$$
e_m(\mathcal K)_X\leq Cm^{-\alpha}[\log_2 n]^{\beta+\alpha}.
$$
Since for $n$ sufficiently large, $2^{-1}\log_2 n<\log_2 m<3\log_2 n$, the proof is completed.
\hfill $\Box$

\begin{remark}
\label{R1}
Similar statement as Lemma \ref{L6} holds if 
$$
d_n(\mathcal{K},n^{an})_X\leq c_0\frac{[\log_2 n]^\beta}{ n^{\alpha}},
$$
where $a>0$ is a positive constant. 
\end{remark}

The results in Lemma \ref{L6} and Remark \ref{R1} hold also for sequences. Namely, using the monotonicity of the quantities involved, the following is true.

\begin{remark}
Let $\mathcal K\subset X$ be a compact subset of a Banach space $X$ with ${\rm rad}(\mathcal{K})<1$.
If  there are constants  $\alpha>0$, $\beta\in \mathbb R$, $\lambda>1$, and $c_0>0$  such that 
 we have
$$
d_n(\mathcal{K},\lambda^n)_X\leq c_0\frac{[\log_2 n]^\beta}{ n^{\alpha}}, \quad n=1,2, \ldots,
$$
or for $a>0$
$$
d_n(\mathcal{K},n^{an})_X\leq c_0\frac{[\log_2 n]^\beta}{ n^{\alpha}},\quad n=1,2, \ldots,
$$
then 
$$
e_n(\mathcal K)_X<C\frac{[\log_2 n]^{\alpha+\beta}}{ n^{\alpha}}, \quad n=1,2, \ldots,
$$
where 
$C$ is a  fixed constant depending only on $\alpha,\beta,c_0$ and $\lambda$ or $a$.
\end{remark}

\section{Applications}
\label{appl}
In this section, we describe how some of the above results can translate to estimates about $m$-term approximation. We follow the 
framework outlined in Theorem 4.1 from \cite{T}.

We assume that we have a system $\mathcal D=\{g_j\}_{j=1}^\infty\subset X$ and
 de la Vallee-Poussin linear operators  $V_k$ associated with the  sequences $n_k$,  $\{(V_k,n_k)\}_{k=1}^\infty$, 
 satisfying the conditions:
\begin{enumerate}
\item There is a constant $A_2>1$ such that 
 $$
V_k(g_j)=\begin{cases}
g_j, &\quad j=1,\ldots,n_k,\\
0, &\quad j>A_2n_k,\\
\alpha_{k,j}g_j, &\quad \hbox{otherwise,}\quad \hbox{where}\quad \alpha_{k,j}\in \mathbb R.
\end{cases}
$$
\item The norms of $V_k$ as operators from $X$ to $X$ are uniformly bounded, i.e. there is a constant $A_3>0$  such that 
 $\|V_k\|_{X\to X}\leq A_3$, $k=1,2,\ldots$.
\end{enumerate}

We denote by $S_{n_k}(f)$ the best approximation to $f\in \mathcal K$ by elements from ${\rm span}\{g_1,\ldots,g_{n_k}\}$, 
$$
E_{n_k}(f,\mathcal D)_X:=\inf_{c_1,\ldots,c_{n_k}}\|f-\sum_{j=1}^{n_k}c_jg_j\|_X=\|f-S_{n_k}(f)\|_X,
$$
and by 
$$
\sigma_m(f,\mathcal D')_X:=\inf_{\{c_j\},\,\Lambda:\,|\Lambda|=m}\|f-\sum_{j\in \Lambda\cap \mathcal D'}c_jg_j\|_X
$$
the best $m$-term approximation of $f$ by a linear combination of $m$ elements from $\mathcal D'$, where
$\mathcal D'$ could be a subset of $\mathcal D$ or $\mathcal D$ itself.
We also define
$$
E_{n_k}(\mathcal K,\mathcal D)_X:=\sup_{f\in \mathcal K}E_{n_k}(f,\mathcal D)_X,
\quad \sigma_m(\mathcal K,\mathcal D')_X:=\sup_{f\in\mathcal K} \sigma_m(f,\mathcal D')_X.
$$
Then the following lemma holds.

\begin{lemma}
\label{TL1}
If the Banach space $X$ admits de la Vallee-Poussin linear operators  $V_k$ that satisfy {\rm (1)-(2)}, with constants $A_2>1$, $A_3>0$, then 
 we have for $1<m<A_2n_k$,
\begin{equation}
\label{klo}
d_m\left(\mathcal K,\left(\frac{A_2bn_k}{m}\right)^m\right)_X\leq (1+2A_3)\max\{E_{n_k}(\mathcal K,\mathcal D)_X,\sigma_{m}(\mathcal K,\mathcal D)_X\},
\end{equation}
where $b>1$ is an absolute constant.
\end{lemma}
\noindent
{\bf Proof:} 
Clearly, we have the inequality
\begin{eqnarray}
\nonumber
\|f-V_k( f)\|_X&\leq& \|f-S_{n_k}(f)\|_X+\|S_{n_k}(f)-V_k(f)\|_X\\
&=&
E_{n_k}(f,\mathcal D)_X+\|V_{k}(S_{n_k}(f)-f)\|_X\leq (1+A_3)E_{n_k}(f, \mathcal D )_X.
\label{to1}
\end{eqnarray}
If we  denote by $\mathcal D_{A_2n_k}:=\{g_1,\ldots,g_{A_2n_k}\}$, then it follows from the properties of $V_k$ that
 for any index set $\Lambda$ with $|\Lambda|=m$ and any coefficients $\{c_j\}_{j=1}^m$, 
$$
\sigma_m(V_k(f),\mathcal D_{A_2n_k})_X\leq \|V_k(f)-V_k(\sum_{j\in \Lambda}c_jg_j)\|_X\leq A_3\|f-\sum_{j\in \Lambda}c_jg_j\|_X,
$$
and therefore
\begin{equation}
\label{to2}
\sigma_m(V_k(f),\mathcal D_{A_2n_k})_X\leq A_3\sigma_m(f,\mathcal D)_X.
\end{equation}
Since
$$
\sigma_m(f,\mathcal D_{A_2n_k})_X\leq\|f-V_k(f)\|_X+\sigma_m(V_k(f),\mathcal D_{A_2n_k})_X,
$$
it follows from (\ref{to1}) and (\ref{to2}) that
\begin{eqnarray*}
\sigma_m(f,\mathcal D_{A_2n_k})_X&\leq& (1+A_3)E_{n_k}(f,\mathcal D)_X+A_3\sigma_m(f,\mathcal D)_X\\
&\leq&(1+2A_3)\max\{E_{n_k}(\mathcal K,\mathcal D)_X,\sigma_m(\mathcal K,\mathcal D)_X\}.
\end{eqnarray*}
Taking a supremum over $f\in \mathcal K$ in the latter inequality gives
\begin{equation}
\label{lo}
\sigma_m(\mathcal K,\mathcal D_{A_2n_k})_X
\leq (1+2A_3)\max\{E_{n_k}(\mathcal K,\mathcal D)_X,\sigma_m(\mathcal K,\mathcal D)_X\}.
\end{equation}
Note that the total number of $m$-dimensional subspaces, $1<m<A_2n_k$, of the linear space ${\rm span}\{g_1,\ldots,g_{A_2n_k}\}$ is
$\binom{A_2n_k}{m}$. Using the Stirling formula, one can show that there is an absolute constant $b>1$ such that 
$$
\binom{A_2n_k}{m}\leq \left(\frac{A_2bn_k}{m}\right)^m.
$$
Then the  definition of nonlinear Kolmogorov width and its monotonicity with respect to $N$ gives
$$
d_m\left(\mathcal K,\left(\frac{A_2bn_k}{m}\right)^m\right)_X\leq 
d_m\left(\mathcal K,\binom{A_2n_k}{m}\right)_X
\leq \sigma_m(\mathcal K,\mathcal D_{A_2n_k})_X.
$$
The latter inequality combined with (\ref{lo}) leads to 
$$
d_m\left(\mathcal K,\left(\frac{A_2bn_k}{m}\right)^m\right)_X\leq (1+2A_3)\max\{E_{n_k}(\mathcal K,\mathcal D)_X,\sigma_{m}(\mathcal K,\mathcal D)_X\},
$$
where $1<m<A_2n_k$, and the proof is completed.
\hfill $\Box$.

We next state a theorem that follows from Lemma \ref{TL1} and our inequalities  for nonlinear Kolmogorov widths in Hilbert spaces.
Note that our theorem  does not require the additional assumptions on the error $E_{n}(\mathcal K,\mathcal D)_H$ that are needed  in Theorem 4.1 from \cite{T} and describes the behavior of the errors in cases not covered by this theorem.
\begin{theorem}
\label{T3neww}
If the Hilbert space $H$ admits de la Vallee-Poussin linear operators  $V_k$ that satisfy {\rm (1)-(2)}, then the following holds:
\begin{itemize}
\item If $e_{n_k}({\mathcal K})_H\geq C \frac{[\log_2n_k]^\beta}{n_k^\alpha},\quad k=1,2,\ldots$, then
there is an absolute constant $C''>0$ such that 
$$
\max\{E_{n_k}(\mathcal K,\mathcal D)_H,\sigma_{m}(\mathcal K,\mathcal D)_H\}\geq C''\frac{[\log_2 m(1+\log_2(A_2n_k/m))]^{\beta-\alpha}}{m^\alpha[1+\log_2(A_2bn_k/m)]^{\alpha}},\quad 
$$
for $1<m<n_k$, $k=1,2,3,\ldots$.
\item If $e_{n_k}({\mathcal K})_H\geq C\frac{1}{[\log_2 n_k]^\alpha},\quad k=1,2,\ldots$, then
there is an absolute constant $C''>0$ such that 
$$
\max\{E_{n_k}(\mathcal K,\mathcal D)_H,\sigma_{m}(\mathcal K,\mathcal D)_H\}\geq C''\frac{1}{[\log_2 m(1+\log_2(A_2bn_k/m))]^\alpha},
$$
for $1<m<n_k$, $k=1,2,3,\ldots$.
\item If $e_{n_k}({\mathcal K})_H\geq C 2^{-cn_k^{\alpha}},\quad k=1,2,\ldots$, then
there are  absolute constants $C''>0$ and $c''>0$ such that 
$$
\max\{E_{n_k}(\mathcal K,\mathcal D)_H,\sigma_{m}(\mathcal K,\mathcal D)_H\}
\geq C''2^{-c''[m(1+\log_2(A_2bn_k/m))]^{\alpha/(1-\alpha)}},
$$
for $1<m<n_k$, $k=1,2,3,\ldots$.
\end{itemize}
\end{theorem}
\noindent
{\bf Proof:}  We use Theorem \ref{T2new}  in the case $N=\left(\frac{A_2bn_k}{m}\right)^m$, Lemma \ref{TL1} and the fact that 
$$
d_m\left(\mathcal K,\left(\frac{A_2bn_k}{m}\right)^m\right)_X\geq 
d_{n_k-1}\left(\mathcal K,\left(\frac{A_2bn_k}{m}\right)^m\right)_X, \quad 1<m<n_k.
$$
Note that we have utilized the fact that the constants  in  Theorem \ref{T2new}  do not depend on $N$. 
\hfill $\Box$

We can derive several corollaries from the above theorem, one of which we state below.
If we take $m=n_k/2$ in Theorem \ref{T3neww}, we obtain the following statement.
\begin{corollary}
If the Hilbert space $H$ admits de la Vallee-Poussin linear operators  $V_k$ that satisfy {\rm (1)-(2)}, then the following holds:
\begin{itemize}
\item If $e_{n_k}({\mathcal K})_H\geq C \frac{[\log_2n_k]^\beta}{n_k^\alpha},\quad k=1,2,\ldots$, then
there is an absolute constant $C''>0$ such that 
$$
\max\{E_{n_k}(\mathcal K,\mathcal D)_H,\sigma_{n_k/2}(\mathcal K,\mathcal D)_H\}\geq C''\frac{[\log_2 n_k]^{\beta-\alpha}}{n_k^\alpha},\quad k=1,2,\ldots.
$$
\item If $e_{n_k}({\mathcal K})_H\geq C\frac{1}{[\log_2 n_k]^\alpha},\quad k=1,2,\ldots$, then
there is an absolute constant $C''>0$ such that 
$$
\max\{E_{n_k}(\mathcal K,\mathcal D)_H,\sigma_{n_k/2}(\mathcal K,\mathcal D)_H\}\geq \frac{C''}{[\log_2 n_k]^\alpha},
\quad k=1,2,\ldots.
$$
\item If $e_{n_k}({\mathcal K})_H\geq C 2^{-cn_k^{\alpha}},\quad k=1,2,\ldots$, then
there are absolute constants $C''>0$, $c''>0$  such that 
$$
\max\{E_{n_k}(\mathcal K,\mathcal D)_H,\sigma_{n_k/2}(\mathcal K,\mathcal D)_H\}
\geq C''2^{-c''n_k^{\alpha/(1-\alpha)}}, \quad k=1,2,\ldots.
$$
\end{itemize}
\end{corollary}

Note that since $A_2>1$, we can take $m=n_k$ in Lemma \ref{TL1}, use the fact that 
$$
E_{n_k}(\mathcal K,\mathcal D)_H\geq \sigma_{n_k}(\mathcal K,\mathcal D)_H,
$$
and obtain from this lemma that if the Banach space $X$ admits de la Vallee-Poussin linear operators satisfying (1)-(2),
then
$$
d_{n_k}\left(\mathcal K,(A_2b)^{n_k}\right)_H\leq (1+2A_3)E_{n_k}(\mathcal K,\mathcal D)_H.
$$
We can now  use Corollary \ref{c1}  and the monotonicity of the nonlinear Kolmogorov width with respect to $N$ to
conclude that 
$$
d_{n_k}\left(\mathcal K,(A_2b)^{n_k+1}\right)_H\leq d_{n_k}\left(\mathcal K,(A_2b)^{n_k}\right)_H,
$$
and derive the following statement.
\begin{corollary}
If the Hilbert space $H$ admits de la Vallee-Poussin linear operators  $V_k$ that satisfy {\rm (1)-(2)}, then the following holds:
\begin{itemize}
\item If $e_{n_k+1}({\mathcal K})_H\geq C \frac{[\log_2n_k]^\beta}{n_k^\alpha},\quad k=1,2,\ldots$, then
there is an absolute constant $C''>0$ such that 
$$
E_{n_k}(\mathcal K,\mathcal D)_H\geq C''\frac{[\log_2 n_k]^{\beta-\alpha}}{n_k^\alpha}, \quad k=1,2,\ldots.
$$
\item If $e_{n_k+1}({\mathcal K})_H\geq C\frac{1}{[\log_2 n_k]^\alpha},\quad k=1,2,\ldots$, then
there is an absolute constant $C''>0$ such that 
$$
E_{n_k}(\mathcal K,\mathcal D)_H\geq \frac{C''}{[\log_2 n_k]^\alpha}, \quad k=1,2,\ldots.
$$
\item If $e_{n_k+1}({\mathcal K})_H\geq C 2^{-cn_k^{\alpha}},\quad k=1,2,\ldots$, then
there are absolute constants $C''>0$, $c''>0$,  such that 
$$
E_{n_k}(\mathcal K,\mathcal D)_H
\geq C''2^{-c''n_k^{\alpha/(1-\alpha)}}, \quad k=1,2,\ldots.
$$
\end{itemize}
\end{corollary}

{\bf Acknowledgment:} This work  was supported by the  NSF Grant  DMS 2134077, Tripods Grant CCF-1934904, and the ONR Contract N00014-20-1-278.

\end{document}